\documentclass[11pt,leqno]{amsart}
\usepackage{amsmath,epsfig}
\usepackage{amssymb}
\numberwithin{equation}{section}
\def\1{\ensuremath{\mathrm{1}\hspace{-.35em} \mathrm{1}}} 

\def\E{\mathbb{E}}

\def\L{\mathbb{L}}
\def\N{\mathbb{N}}
\def\P{\mathbb{P}}
\def\R{\mathbb{R}}
\def\Z{\mathbb{Z}}
\def\cov{\mathop{\rm Cov}\nolimits}
\def\Lip{\mathop{\rm Lip}\nolimits}

\renewcommand{\le}{\ensuremath{\leqslant}}
\renewcommand{\ge}{\ensuremath{\geqslant}}
\renewcommand{\epsilon}{\varepsilon}

\newcommand{\refeqn}[1]{eqn$.$~(\ref{#1})}

\def\ie{{\it i.e{$.$}}}

\newcommand{\introo}[2]{{\left]{#1,\,#2\,}\right[\kern1pt}}

\newcommand{\intrfo}[2]{{\left[{#1,\,#2}\right[\kern1pt}}

\newtheorem{theo}{Theorem}[section]
\newtheorem{Lem}{Lemma}[section]
\newtheorem{Prop}{Proposition}[section]

\newtheorem{Cor}{Corollary}[section]
\newtheorem{Rem}{Remark}[section]

\begin{document}

\title[Weakly dependent chains with infinite memory]
{Weakly dependent chains with infinite memory}
\author[P. Doukhan]{Paul Doukhan$^{(1),(2)}$}
\address[1]{LS-CREST,
Laboratoire de Statistique,
Timbre J340,
3, avenue Pierre Larousse,
92240 MALAKOFF, FRANCE}

\email{doukhan@ensae.fr}

\author[O. Wintenberger]{Olivier Wintenberger$^{(2)}$}
\address[2]{SAMOS-MATISSE
(Statistique Appliqu\'ee et MOd\'elisation Stochastique)
Centre d'\'Economie de la Sorbonne
Universit\'e Paris 1 - Panth\'eon-Sorbonne, CNRS
90, rue de Tolbiac
75634-PARIS CEDEX 13, FRANCE}

\email{olivier.wintenberger@univ-paris1.fr}


\medskip
\begin{abstract}
We prove the existence of a weakly dependent strictly stationary
solution of the equation $
X_t=F(X_{t-1},X_{t-2},X_{t-3},\ldots;\xi_t)$ called {\em chain with
infinite memory}. Here the {\em innovations} $\xi_t$ constitute an
independent and identically distributed sequence of random variables. The function $F$ takes values in
some Banach space and satisfies a Lipschitz-type condition. We also
study the interplay between the existence of moments and the rate of
decay of the Lipschitz coefficients of the function $F$. With the
help of the weak dependence properties, we derive Strong Laws of
Large Number, a Central Limit Theorem and a Strong Invariance
Principle.
\end{abstract}

\maketitle
\section{Introduction}
%
\noindent Statistical inferences heavily rely on the underlying model. 
Processes may have different representations. Thus  they belong to
different classes of models.
In this paper,
we introduce a {\em chain with infinite memory} as the stationary
solution of the equation
\begin{equation}\label{chaineinf}
 X_t=F(X_{t-1},X_{t-2},X_{t-3},\ldots;\xi_t), \qquad
\mbox{ {\it a.s.} \qquad for }\quad t\in\Z,
\end{equation}
where $F$ takes value in a Banach space. For details, see section \ref{sec::model}.
The dynamical behavior described by \eqref{chaineinf} corresponds to a large variety of times series models. Those models can be seen as natural extensions, either of linear models or of Markov models. In the sequel,  the {\em innovations} $\xi_t$
constitute an independent and identically distributed (iid) sequence. Various representations use such
innovations. For instance, the case of causal Bernoulli shifts
$X_t=H(\xi_t,\xi_{t-1},\ldots)$ was studied by Wu
\cite{Wu2005}. But several Bernoulli shifts, such as Volterra series,
may not fit the parsimony criterion and the function $H$ may be non-explicit. This is a drawback for statistical inferences in that
context. Markov models are preferred in various applications e.g. in
finance, hydrodynamics, physics, see \cite{DOKU93,LALI80}.
Kallenberg \cite{KALL97} stresses the fact that all the $p$-Mar\-kov processes
are solutions of equations of the type:
\begin{equation} \label{markov}
 X_t=F(X_{t-1},\ldots,X_{t-p};\xi_t).
\end{equation}
Bougerol \cite{BOUG93} gave conditions of Lyapunov type for the
existence of a stationary solution to {\em Stochastic Recurrence Equations} (SRE), which are particular cases of \eqref{markov}.\\

\noindent Other approaches than \eqref{chaineinf} to modeling processes which do not satisfy the
Markov property already exist; the \emph{Random Systems with
Complete Connections} (RSCC), see \cite{IOTH69}, and
the \emph{Variable Length Markov Chains} (VLMC), see
\cite{BUWY99}. Such models are widely used in the fields of
particle systems or in DNA data analysis.  These processes are
defined through their conditional distributions. Their existence
relies on assumptions on the conditional expectations, following the
work of Dobrushin \cite{DOBR70}. Notice that Berbee \cite{BERB87} obtained
another existence condition for the cases where the state space is
discrete, see also \cite{COFE02,FEMA04}.\\

\noindent Dobrushin's condition implies
mixing, see \cite{DOUK94,IOTH69}. Mixing coefficients are useful to derive asymptotic theorems for various functionals of a stationary sequence, see
Rio \cite{Ri}. However, major asymptotic results still hold under so-called weak
dependence conditions, see section \ref{wd}, \cite{DEPR04,DOLO99}
and the recent monograph by \cite{DEDO05}. The Central Limit Theorem
(CLT) in Dedecker and Doukhan \cite{DEDO03} holds if the $x^2\ln (1+x)$th moments of
$X_0$ are finite and if the process is weakly dependent with
geometric decay of the coefficients. Because weak dependence is less
restrictive than mixing (see Andrews \cite{ANDR84} for an example) this
result extends the CLT for mixing sequences due to Rio \cite{Ri}. The conditions for the CLT are
expressed in terms of Orlicz functions that balance the moments of some order and the weak dependence conditions.\\

\noindent The existence of a stationary solution to \eqref{chaineinf} is
proved in section \ref{exist} under a specific Lipschitz type
assumption on $F$, see \eqref{lipcondf} below. Approximation by
suitable Markov processes is the main tool for the proofs given in section
\ref{proofs}. This existence condition also yields finiteness of
moments of some order in terms of Orlicz functions. We get bounds
for the weak dependence coefficients of the solution to \eqref{chaineinf}. We use these
bounds to derive sufficient conditions on $F$ in term of Orlicz functions and in turn to prove a
Strong Law of Large Numbers (SLLN), a CLT and a Strong Invariance
Principle (SIP), see section \ref{asres}. We discuss the generality of our model in section \ref{secexample} comparing it with existing ones. But to begin with, we
introduce some notation and we define useful tools such as weak
dependence and Orlicz spaces.
\section{Preliminaries}\label{prel}
\subsection{Notation}
In the sequel, the iid innovations $\xi_t$ for $t\in\Z$  take
values in a measurable space $(E',\mathcal{A}')$. Let $\|\cdot\|$ denote the norm of a Banach space $E$. The space $E^{(\infty)}$ is the subset of $E^\N$ of finitely-non-zero sequences $(x_k)_{k>0}$ such that there exists $N> 0$ with $x_k=0$ for $k>N$. Let $E$ be endowed with its Borel $\sigma-$algebra $\mathcal{A}$, then $E^{(\infty)}$ is considered together with its product $\sigma-$algebra $\mathcal{A}^{\otimes \N}$. The function $F$ in \eqref{chaineinf} is assumed to be  a measurable function from $E^{(\infty)}\times E'$ with value in $E$.  Moreover $\|\cdot\|_m$ denotes the usual $\L^m$-norm, \ie, $\|X\|_m^m=\E \|X\|^m$ for $m\ge 1$ for every $E$-valued random variable $X$.  For
$h:E\to \R$, we denote $\|h\|_\infty=\sup_{x\in E}|h(x)|$ and
$$
\Lip(h)=\sup_{x\ne y} \frac{|h(x)-h(y)|}{\|x-y\|}.
$$
The space $\Lambda_1\big(E\big)$  is the set of functions $h:E\to\R$
such that $\Lip(h)\le1$.

\subsection{Weak dependence}\label{wd}

An appropriate notion of weak dependence for the model \eqref{chaineinf} was introduced in \cite{DEPR04}.
It is based on the concept of the coefficient $\tau$ defined below.
Let $(\Omega, \mathcal{C}, \P)$ be a probability space,
$\mathcal{M}$ a $\sigma$-subalgebra of $\mathcal{C}$ and $X$ a
random variable with values in $E$. Assume that $\|X\|_1<\infty$ and define the coefficient $\tau$ as
\begin{equation*}
\tau(\mathcal{M},X)=\\
\left\|\sup\left\{\left|\int f(x)\P_{X|\mathcal{M}}(dx)-\int
f(x)\P_{X}(dx)\right|\mbox{ with } ~f\in\Lambda_1\big(E\big)\right\}\right\|_1.
\end{equation*}
\label{coupling}\noindent An easy way to bound this coefficient is based on a coupling
argument: $$\tau(\mathcal
M,X)\le\|X-Y\|_1$$ for any $Y$ with the same distribution as $X$ and independent
 of $\mathcal{M}$, see \cite{DEPR04}.
Moreover, if the probability space $(\Omega, \mathcal{C}, \P)$ is rich enough (we always assume so in the sequel)  there exists an $X^\ast$ such that $\tau(\mathcal M,X)=\|X-X^\ast\|_1$.
Using the definition of $\tau$, the dependence between the past of the sequence $(X_t)_{t\in\Z}$ and its future $k$-tuples may be assessed: Consider the norm
$\|x-y\|=\|x_1-y_1\|+\cdots+\|x_k-y_k\|$ on $E^k$,
set $\mathcal{M}_p=\sigma(X_t,t\le p)$ and define
\begin{eqnarray*}
\tau_k(r)&=&\max_{1\le l\le k}
\frac1l\sup\Big\{\tau(\mathcal{M}_p,(X_{j_1},\ldots,X_{j_l}))\mbox{ with }p+r\le
j_1<\cdots <j_l\Big\},\\
\tau_\infty(r)&=&\sup_{k>0}\tau_k(r).
\end{eqnarray*}
For the sake of simplicity, $\tau_\infty(r)$ is denoted by $\tau(r)$.
Finally, the time series $(X_t)_{t\in\Z}$ is {\it $\tau$-weakly dependent} when its coefficients $\tau(r)$ tend to $0$ as $r$ tends to infinity.

\subsection{Orlicz spaces}
Orlicz spaces are convenient generalizations of the classical $\L^m$-spaces, we refer to \cite{KR} for the introduction and properties of such spaces.
Let $\Phi$ be an Orlicz function, \ie, defined on $\R^+$, convex, increasing and satisfying $\Phi(0)=0$. For any random variable $X$ with values in $E$, the norm $\|X\|_{\Phi}$ is defined by the equation
$$\|X\|_{\Phi}=\inf\left\{u>0\mbox{ with }\E\Phi\left(\frac{\|X\|}{u}\right)\le 1\right\}.$$
The Orlicz space $\L^{\Phi}$ is given by
$$
\L^{\Phi}=\big\{E\mbox{-valued random variables }X\mbox{ such that } \|X\|_\Phi<\infty\big\}.
$$
It is a Banach space equipped with the norm $\|\cdot\|_\Phi$. For $m\ge 1$ and $\Phi(x)=x^m$, notice that $\L^\Phi$ is the usual $\L^m$-space. We restrict ourselves to Orlicz functions $\Phi$ satisfying the condition:
\begin{equation}\label{submult}
\mbox{For all }x,y\in {\R^+},\qquad \Phi(xy)\le \Phi(x)\Phi(y).
\end{equation}
This class of Orlicz functions is sufficiently large. For instance, the functions $\Phi(x)=x^m$ and $\Phi(x)=x^m(1+\ln(1+x))^{m'}$ satisfy \eqref{submult} for each $m\ge 1,m'\ge 0$. Moreover, if $\phi$ is any  Orlicz function satisfying the $\Delta_2$-condition (there exists $k>0$ such that $\phi(2x)\le k\phi(x)$) then $\Phi(x)=\sup_{u>0}\phi(xu)/\phi(u)$ is an  Orlicz function satisfying \eqref{submult}. Various examples of Orlicz functions satisfying the $\Delta_2$-condition are given in \cite{KR}.\\

\noindent Later, in theorem \ref{clt} we will need some transformations of Orlicz functions. Given such a function $\Phi$, we define for $q> 1$,
\begin{equation}\label{tildephi}
  \widetilde\Phi_q(x)=\sup_{y>0}\{(xy)^{q-1}-\Phi(y)/y\}.
\end{equation}
The transformations $\widetilde\Phi_q(x)$ have simple bounds for certain choices of $\Phi$, see lemma \ref{calctr} for details. In particular,
if $\Phi(x)=x^m$ for $m>q>1$, then $  \widetilde\Phi_q(x)\le x^{(m-1)(q-1)/(m-q)}$.
Another useful example is the one of $\Phi(x)=x^q(1+\ln(1+x))^{(1+b)(q-1)}$ and $  \widetilde\Phi_q(x)\le \exp((q-1)x^{1/(1+b)})x^{q-1}$ for any $q>1$ and $b\ge 0$.

\section{The results}\label{results}
\subsection{Assumptions}\label{sec::model}
The existence of a solution to \eqref{chaineinf} will be proved under a Lipschitz-type condition. We express it in terms of some Orlicz functions in order to be able to work with moments more general than power moments, see theorem \ref{thchaineinf}. These moments will be needed to establish the asymptotic results of theorem~\ref{clt}.\\

\noindent Assume there exists an Orlicz function $\Phi$ such that for all $x$, $y$ in $E^{(\infty)}$:
\begin{equation}\label{lipcondf}
\left\|F(x;\xi_0)-F(y;\xi_0)\right\|_\Phi
\le\sum_{j=1}^\infty a_j\|x_j-y_j\|,
\end{equation}
where $(a_j)_{j\ge1}$ is a sequence of nonnegative real numbers such that
\begin{eqnarray}
\label{propinf} a~&=&\sum_{j=1}^\infty a_j\,<\,1\mbox{ and }
\\
\label{propmom}\mu_\Phi&=&\left\|F(0,0,\ldots;\xi_0)\right\|_\Phi\,<\, \infty.
\end{eqnarray}
The Lipschitz property of $F$ and the moment assumption \eqref{propmom} induce that $\|F(c;\xi_0)\|_\Phi<\infty$ for any constant $c\in E^{(\infty)}$. We choose $c=(0,0,\ldots)$ in condition \eqref{propmom} for convenience.

\subsection{Existence, moments and weak dependence}\label{exist}
The following theorem settles the existence of a solution to \eqref{chaineinf}. It also states that the $\Phi$th moment of this solution is finite.
\begin{theo}\label{thchaineinf}
Assume that conditions \eqref{propinf} and \eqref{propmom} hold for some Orlicz function
$\Phi$ satisfying \eqref{submult}. Then there exists a $\tau$-weakly dependent stationary solution $(X_t)_{t\in\Z}$ of \eqref{chaineinf} such that $\|X_0\|_\Phi<\infty$ and
$$\tau(r)\le 2\frac{\mu_1}{1-a}\inf_{1\le p\le r}\left(a^{r/p}+\frac{1}{1-a}\sum_{k=p+1}^\infty  a_{k}\right)\to 0\quad\mbox{ as }\quad r\to\infty.$$
\end{theo}
\noindent The proof of the existence of a solution to \eqref{chaineinf} is given in section \ref{sec::chaineinf} expressing it as the limit of the $p$-Markov processes defined in \eqref{markov}. The weak dependence properties are proved in section \ref{wd::chaineinf}.

\begin{Rem}\label{uniq}\normalfont  We also prove in section \ref{proofs} that there exists some measurable function $H$ such that $X_t=H(\xi_t,\xi_{t-1},\ldots)$. This means that the process $(X_t)_{t\in\Z}$ can be represented as a causal Bernoulli shift. For those processes, conditions \eqref{propinf} and \eqref{propmom} together imply the Dobrushin uniqueness condition, see \cite{DOBR70}. Thus $(X_t)_{t\in\Z}$ is the unique causal Bernoulli shift solution to \eqref{chaineinf}. Moreover, as a causal Bernoulli shift, the solution $(X_t)_{t\in\Z}$ is automatically an ergodic process. Under the conditions of theorem \ref{thchaineinf}, the solution to \eqref{chaineinf} has finite $\Phi$th moment. From lemma \ref{phix}, $(X_t)_{t\in\Z}$ has also finite first order moments. The ergodic theorem yields the SLLN for any {\em chain with infinite memory} under the assumptions of theorem \ref{thchaineinf}.
\end{Rem}

\begin{Cor}\label{The finite memory case}
Under the assumptions of theorem \ref{thchaineinf}, there exists a $\tau$-weakly dependent stationary solution $(X_t)_{t\in\Z}$ to \eqref{markov}
such that $\|X_0\|_\Phi<\infty$ and
$\displaystyle\tau(r)\le 2\mu_1(1-a)^{-1}~a^{r/p}$ for $r\ge p$.
\end{Cor}
\noindent Dedecker and Prieur \cite{DEPR04} proved the existence of a solution to \eqref{markov}. They stated that there exists $0<\rho<1$ and $C>0$ such that $\tau(r)\le C\rho^r$. Applying corollary \ref{The finite memory case}, we get the bound $\rho\le a^{1/p}$.
The bounds of the weak dependence coefficients in theorem \ref{thchaineinf} come from an approximation with Markov chains of order $p$ and from the result of corollary \ref{The finite memory case}.\\

\noindent In theorem \ref{thchaineinf}, the $\tau$-weak dependence property is linked to the choice of the parameter $p$ and then to the rate of decay of the Lipschitz coefficients $a_j$. For example, if $a_j\le ce^{-\beta j}$, we choose $p$ as the largest integer smaller than $\sqrt{-\ln(a)
r/\beta}$ to derive the bound
$\tau(r)\le C e^{-\sqrt{-\ln(a)\beta r}}$ for some suitable
constant $C>0$. If $a_j\le cj^{-\beta}$, we choose the largest integer $p$ such that $p\ln p(1-\beta)/\ln a\le r$. Then there exists $C>0$ such that $\tau(r)\le C p^{1-\beta}$. Notice that $\ln r$ is smaller than $\ln p+ \ln \ln p$ up to a constant and that $\ln r/r$ is proportional to  $1/p(1+\ln\ln p/\ln p)$ and then equivalent to $1/p$ as $p$ tends to infinity with $r$. From these equivalences, we achieve thus there exists $C>0$ such that $\tau(r)\le C\left(\ln r/r\right)^{\beta-1}$. \\

\noindent A similar result as the one of theorem \ref{thchaineinf} was obtained for discrete state space models (as RSCC) in \cite{IOGR90}. They gave bounds for the mixing coefficients under conditions on the marginal
distributions of the {\em innovations}. The bound
in \cite{IOGR90}, theorem 2.1.5 on page 42,  is similar to the one for $\tau(r)$ in theorem \ref{thchaineinf}. In a sense we extend their result: Here the {\em innovations} are not supposed to be absolutely continuous and our approach can be applied to discrete state space processes as well, see the example of the {\em Galton-Watson process with immigration} in section \ref{secexample}.\\

\noindent Bougerol gives in \cite{BOUG93} a recursive approximation of the stationary measure in the Markovian case. In proposition \ref{simulation} below we generalize this result to the infinite memory case. Let $\phi_k:E^{k-1}\times E'\to E$ be the random function defined as $x\mapsto F(x,c;\xi_k)$, for each $k\ge2$ and some fixed sequence $c=(c_1,c_2,\ldots)\in E^{(\infty)}$. Write $\widetilde X_1=\phi_1=\phi(c;\xi_0) $ and define recursively
$$
\widetilde X_n=\phi_n(\widetilde X_{n-1},\ldots, \widetilde X_1).
$$
\begin{Prop}\label{simulation}
Assume that conditions \eqref{propinf} and \eqref{propmom} hold for $\Phi$  satisfying \eqref{submult}.
If $(X_t)_{t\in\Z}$ is the solution to \eqref{chaineinf} then
$$ \|\widetilde X_r-X_{r}\|_\Phi\le \frac{\|X_0\|_\Phi+\overline c}{1-a}\inf_{1\le p\le r}\left(a^{r/p}+\frac{1}{1-a}\sum_{k=p+1}^\infty  a_{k}\right)\to 0\quad\mbox{ as }\quad r\to\infty,$$
where $\overline c$ is a constant such that $\|c_i\|\le \overline c$ for all $i\ge1$.
\end{Prop}
\noindent The proof of this proposition is given in section \ref{pfsimu}.

\subsection{Asymptotic results}\label{asres}
In this section, $E=\R$. We give an appropriate condition on $F$ (see \eqref{Dp} below) that leads to versions of the results of Dedecker and Doukhan \cite{DEDO03} and Dedecker and Prieur \cite{DEPR04} obtained under weak dependence.
\begin{theo}\label{clt}\label{flil}
Assume that conditions \eqref{propinf} and \eqref{propmom} hold for some Orlicz function $\Phi$  satisfying \eqref{submult} and assume there exists $c_0>0$ such that
\begin{subequations}\label{Dp}
\begin{align}
& \displaystyle\sum_{k\ge1}a_k   \widetilde\Phi_q
(c_0k)<\infty \quad\mbox{ if there exist }p\ge1\mbox{ such that }\sum_{j>p}a_j=0,\label{Dp1}\\
&\displaystyle \sum_{k\ge1}a_k   \widetilde\Phi_q
\left(-c_0k\ln\Big(\sum_{j\ge k}a_j\Big)\right)<\infty\label{Dp2}
\quad\mbox{ otherwise,}
 \end{align}
\end{subequations}
where $\widetilde \Phi_q$ is defined in \eqref{tildephi}. The following relations hold:
\begin{description}
\item[SLLN] If $q\in]1,2[$ then \qquad$\displaystyle n^{-1/q}\sum_{i=1}^{n}(X_i-\E X_0) \to_{n\to\infty}0,$\qquad a.s.
\item[CLT] If $q= 2$, then\qquad
$\displaystyle
\frac1{\sqrt n}\sum_{i=1}^{[nt]}(X_i-\E X_0)\stackrel{D[0,1]}{\longrightarrow}\sigma W(t)$\qquad as $n\to\infty$\\ where  $\displaystyle\sigma^2=\sum_{i=-\infty}^\infty\cov(X_0,X_i)$ is finite and $W(t)$ is the standard Wiener process.
\item[SIP] If $q=2$ and if the underlying probability space is rich enough then there exist independent $\mathcal{N}(0,\sigma^2)$-distributed random variables $(Y_i)_{i\ge 1}$ such that\\
$$
\sum_{i=1}^n(X_i-Y_i)=o(\sqrt{n\ln \ln n}) \qquad  \mbox{a.s.}
$$
\end{description}
\end{theo}
\noindent The proof of this theorem is given in section
\ref{pfclt}.\\

\noindent Note that $x^2\ln(1+ x)$th moments are necessary to get
the CLT for weakly dependent processes. See \cite{dmr} for an
example of processes, solutions of \eqref{markov} for $p=1$, that do
not satisfy the CLT under conditions \eqref{propinf} and
\eqref{propmom} for $\Phi(x)=x^2$. Note also that approximations by
martingale difference as in \cite{Peligrad2005} or projective
criterion as in \cite{Merlevede2006} give the CLT under weaker
assumptions for some
of the examples treated in section \ref{secexample}.\\

\noindent Condition \eqref{Dp1} is relevant for the Markov solution
$(X_t)_{t\in\Z}$  to \eqref{markov}, {\it i.e.}, when
$\sum_{j>p}a_j=0$. For the other cases, we rewrite assumption
\eqref{Dp2} for various rates of decay of the Lipschitz coefficients
$a_j$. Let $a,b,c$ be some positive real numbers then
\begin{align}
\label{Dp11}\tag{\ref{Dp1}'} &\mbox{If }a_k\le ck^{-a}, &\sum_{k\ge1}a_k
\widetilde\Phi_q\left(c_0k\ln k\right)<\infty
\qquad \mbox{ for some  }c_0>0.\\
\label{Dp12}\tag{\ref{Dp1}''}&\mbox{If }a_k \le c \exp(-ak^b),
&\sum_{k\ge1}a_k   \widetilde\Phi_q\left(c_0k^{1+b}\right)<\infty
\qquad\mbox{ for some }c_0>0.
\end{align}
For instance, condition \eqref{Dp11} holds if $\Phi(x)=x^m$ for $m>q$ and
 $\displaystyle a>1+(q-1)(m-1)(m-q)^{-1}$.
Condition \eqref{Dp12} holds for $\Phi(x)=x^q(1+\ln(1+x))^{(1+b)(q-1)}$. Applying theorem \ref{clt}, the CLT and the SIP hold for sub-geometric rates of decay of the Lipschitz coefficients as in \eqref{Dp12} under a moment condition of the order $x^2(1+\ln(1+x))^{1+b}$.

\section{Examples}
\label{secexample}
\noindent In this section, we present some examples with $E=\R^d$ and $d\ge1$. We consider the finite memory case as well as an
 infinite memory extension of  {\em Stochastic Recurrence Equations} (SRE). In particular, we consider the example of the {\em Galton-Watson process with immigration} which satisfies the conditions of our results, but it is not a SRE  in the sense of \cite{BOUG93}.
\subsection{Markov models}

\begin{itemize}
\item[] \textbf{SRE.}\
We consider an iid process $(\phi_t)_{t\in\Z}$ of random Lipschitz maps with $\|\phi_t(x_1)-\phi_t(y_1)\|\le L(\phi)\|x_1-y_1\|$ a.s. for all $x_1,y_1 \in E$ and $t\in\Z$. Moreover let $\phi_t(x)$ be measurable for every fixed $x\in E$ and $t\in\Z$. If a stochastic process $(X_t)_{t\in\Z}$ with values in $E$ satisfies the equation
$$
X_{t+1}=\phi_t(X_t)~~~~\mbox{ a.s., for all }t\in\Z,
$$
we say that $(X_t)_{t\in\Z}$ obeys the SRE associated with $(\phi_t)_{t\in\Z}$. We write this equation as in \eqref{markov} setting $\xi_t=\phi_t$ for $t\in\Z$, and $F(x,z)=z(x)$ for $x\in E$ and $z\in E'$, the space of Lipschitz random functions. In this case, conditions \eqref{propinf} and \eqref{propmom} become
$$
\| L(\phi)\|_\Phi<1\quad\mbox{ and }\quad\|\phi_0(0)\|_\Phi<\infty.
$$
Weaker conditions related to a Lyapunov exponent for the existence of an {\it a.s.} solution to a SRE are obtained in \cite{BOUG93}. However, that result does not yield the existence of moments nor asymptotic results as those in theorem \ref{clt}.
We also mention the survey article by \cite{DIFR99} for an overview and other nice application of SREs.
\item[] \textbf{Nonlinear autoregressive models.}\ Here we consider a solution to \eqref{chaineinf}, where  ${E'}=E$  and $F$ admits the representation
$$F(x_{1},\ldots,x_{p};s)=R(x_{1},\ldots,x_{p})+s.$$
Condition \eqref{propinf} becomes
$$
\|R(y_{1},\ldots,y_{p})-R(x_{1},\ldots,x_{p})\|\le
\sum_{j=1}^{p}a_j\|x_j-y_j\|\mbox{ with }\sum_{j=1}^{p}a_j<1,
$$
and condition \eqref{propmom} coincides with
$\|\xi_0\|_\Phi<\infty$. Results similar to those in theorem \ref{clt} are obtained by different methods in   \cite{DUFL97}.
\item[] \textbf{Galton-Watson processes with immigration.}
If $E=\R$, a {\em Galton-Watson process with immigration} is given as a stationary solution of the equation
\begin{equation}\label{gwi}
X_t=\begin{cases}\displaystyle
\sum_{i=1}^{X_{t-1}} \zeta_{t,i}+\zeta_t,&\quad{\mbox{ if
}}\quad X_{t-1}>0,\\
\displaystyle 0 &\quad{\mbox{ if
}}\quad X_{t-1}=0.\end{cases}
\end{equation}
Here $(\zeta_{t,i})_{t\in\Z,i>0}$, $(\zeta_{t})_{t\in\Z}$ are independent iid families of integer-valued random variables and $E'=\N^{\N}$ is equipped with the product measure. We can write $X_t=F(X_{t-1},\xi_t)$ with
$F(x,(u_i)_{i\ge0})=u_0+\sum_{i=1}^xu_i$ if $x>0$ and $F(0,(u_i)_{i\ge0})=0$ for any $(u_i)_{i\ge0}$.
If $y_1>x_1>0$ then $F(x,(u_i)_{i\ge0})-F(x,(u_i)_{i\ge0})=\sum_{i=x_1}^{y_1}u_i$ thus
$$\|F(x_1,\xi_0)- F(y_1,\xi_0)\|_\Phi=\big\|\sum_{i=x_1}^{y_1}\zeta_{0,i}\big\|_\Phi=|y_1-x_1|\|\zeta_{0,0}\|_\Phi.$$ Assumptions \eqref{propinf} and \eqref{propmom} hold as soon as $\|\zeta_{0,0}\|_\Phi<1$. This model is not a SRE if $\zeta_{0,0}$ is not finitely supported, thus we are not under the conditions of \cite{BOUG93}. Other non-SRE examples which can be treated by our approach are given in \cite{Latour1997}.
\end{itemize}

\subsection{SRE with infinite memory}
Infinite memory extensions of classical SRE are solutions of the equation
$$
\begin{cases}\displaystyle &X_t=\phi_t(X_{t-1},X_{t-2},\ldots) \qquad a.s.,\\
\displaystyle &\|\phi_t(x)-\phi_t(y)\|\le \sum_{i=1}^\infty L_i(\phi)\|x_i-y_i\|,\qquad a.s.
\end{cases}$$ for all $x=(x_i)_{i\ge1},y=(y_i)_{i\ge1} \in E^{(\infty)}$.
Here $(\phi_t)_{t\in\Z}$ is an iid process of random Lipschitz maps.
If $\sum_{i\ge 1}\|L_i(\phi)\|_\Phi<1$ then conditions \eqref{propinf} and \eqref{propmom} are satisfied. Some examples with this representation follow.
\begin{itemize}
\item[] \textbf{Non Linear ARCH$(\infty)$  models.}
\normalfont Here $(X_t)_{t\in\Z}$ is the stationary solution of the equation
$$
X_t=\xi_t\left( \alpha+\sum_{j=1}^\infty  \alpha_j(X_{t-j})\right),
$$
where $\xi_t$ is a $d\times k$ matrix, ${E'}=\mathcal{M}_{k,d}(\R)$,
$\alpha\in\R^k$ and $\alpha_j:E\to\R^k$ are Lipschitz functions. The
LARCH$(\infty)$ model of \cite{GILE04,DOTE06} corresponds to the
special case of linear functions $\alpha_j(x)=c_jx$ with $k\times d$ matrices $c_j.$ Assumptions \eqref{propinf} and
\eqref{propmom} hold as soon as $\|\xi_0\|_\Phi\sum_{j\ge 1}  \Lip \alpha_j<1$ and $\sum_{j\ge 1}
\alpha_j(0)<\infty$.
\item[] \textbf{Models with linear input.}
\normalfont
Let $f:\R^k\times {E'}\to E$ be measurable and satisfy
 $\|f(t,\xi_0)-f(s,\xi_0)\|_\Phi\le L\;\|t-s\|$ for some finite constant $L>0$. We consider
 $$
X_t=f(A_t,\xi_t),\qquad A_t=\sum_{j=1}^\infty c_jX_{t-j},
 $$
where $c_j$ are $k\times d$ matrices. Relations
\eqref{propinf} and \eqref{propmom} hold if $L \sum_{j\ge 1} \|c_j\|<1$ and $f(0,\xi_0)\|_\Phi<\infty$. These models
are used in statistical mechanics, see \cite{KACM59}.

\item[] \textbf{Affine models.}
\normalfont
Let us consider the special case of chains with infinite memory that
can be written in a bilinear form
\begin{equation}\label{modaf}
X_t=M_t \xi_t+f_t,
\end{equation}
where $M_t=M(X_{t-1},X_{t-2},\ldots)$ and $f_t=f(X_{t-1},X_{t-2},\ldots)$ are both
Lipschitz functions of the past values
$X_{t-1},X_{t-2},X_{t-3},\ldots$. Applying theorem \ref{thchaineinf} under the condition
$$
\|\xi_0\|_\Phi\sum_{i=1}^\infty\Lip M_i+\sum_{i=1}^\infty\Lip f_i<1,
$$ there
exists a weakly dependent solution to \eqref{modaf}. This class
contains various time series models (such as ARCH, GARCH, ARMA,
ARMA-GARCH, etc.). In the appendix we prove  the existence of
the joint densities of the solution to \eqref{modaf}. This result
and the weak dependence properties obtained in theorem
\ref{thchaineinf} are needed for achieving optimal rates of
convergence of nonparametric estimators, see \cite{RAWI06}.
\end{itemize}
\section{Proofs of the main results}\label{proofs}
\noindent After some preliminaries in section \ref{preli}, in section
\ref{approx} we construct  a solution of the Markov model \eqref{markov}. We
use it to approximate the solution to \eqref{chaineinf}.  The
existence of a solution to \eqref{chaineinf}, presented in theorem \ref{thchaineinf}, is obtained as $p\to\infty$
in section \ref{sec::chaineinf}. Its weak dependence properties
are derived by coupling techniques in section \ref{wd::chaineinf}. Using weak
dependence results of \cite{DEDO03,DEPR04}, we prove theorem
\ref{clt} in section \ref{pfclt}. Finally, we derive the proof of proposition
\ref{simulation} in section \ref{pfsimu}.

\subsection{Preliminaries}\label{preli}
We first present four useful lemmas. The first one aims at bounding the transformations $\tilde \Phi_q$ for $q>1$, the other ones
are used in the proof of the existence of a solution of
\eqref{chaineinf}.
\begin{Lem}\label{calctr}
Assume $L$ is an increasing non-negative function on $[0,\infty]$ and write $L^{-1}$ for the generalized inverse of $L$, {\it i.e.}, $L^{-1}(x)=\inf\{y>0~\mbox{ with }~L(y)\ge x\}$. If
$\Phi(x)=x^qL(x)$, $x\ge0$, for some $q>1$ then
$$
\widetilde\Phi_q(x)\le \left(xL^{-1}(x^{q-1})\right)^{q-1}\mbox{ for all }x\ge0,
$$
\end{Lem}
\begin{proof} From \eqref{tildephi}, we have $\widetilde\Phi_q(x)=\sup_{y>0}\{y^{q-1}(x^{q-1}-L(y))\}$. We restrict ourselves to $y\le L^{-1}(x^{q-1})$ otherwise $y^{q-1}(x^{q-1}-L(y))\le 0$. Now notice that the first term of the product $y^{q-1}(x^{q-1}-L(y)$ is increasing and the second term always remains smaller than $x^{q-1}$. This proves the lemma.
\end{proof}
\begin{Lem}\label{preliminary2}
Assume that the Orlicz function $\Phi$ satisfies \eqref{submult}. Let $\xi$ and $\zeta$ be independent random variables, $z$ a measurable function and $Z=z(\xi,\zeta) $. We write $\E_\xi$ for the expectation with respect to the distribution of $\xi$. Define
\begin{equation}\label{normxi}
\|z(\xi,\zeta)\|_{\Phi,\xi}=\inf\left\{u>0\mbox{ with }~~\E_\xi\Phi(\|z(s,\zeta)\|/u)\le 1\right\}.
\end{equation}
Then  $\|Z\|_\Phi\le \left\|\|Z\|_{\Phi,\xi}\right\|_\Phi.$
\end{Lem}
\begin{proof}
One needs to prove that $\E\Phi(Z/\|\|Z\|_{\Phi,\xi}\|_\Phi)\le 1$:
\begin{eqnarray*}
\E\Phi\left(\frac{Z}{\|\|Z\|_{\Phi,\xi}\|_\Phi}\right)\le\E\Phi\left(\frac{Z}{\|Z\|_{\Phi,\xi}}\frac{\|Z\|_{\Phi,\xi}}{\|\|Z\|_{\Phi,\xi}\|_\Phi}\right)\le\E\left[\Phi\left(\frac{Z}{\|Z\|_{\Phi,\xi}}\right)\Phi\left(\frac{\|Z\|_{\Phi,\xi}}{\|\|Z\|_{\Phi,\xi}\|_\Phi}\right)\right].
\end{eqnarray*}
The last inequality follows from \eqref{submult}. By independence of $\xi$ and $\zeta$ and by \eqref{normxi}
\begin{eqnarray*}
\E\Phi\left(\frac{Z}{\|\|Z\|_{\Phi,\xi}\|_\Phi}\right)\le\E\left[\Phi\left(\frac{\|Z\|_{\Phi,\xi}}{\|\|Z\|_{\Phi,\xi}\|_\Phi}\right)\E_{\xi}\Phi\left(\frac{Z}{\|Z\|_{\Phi,\xi}}\right)\right]\le\E\Phi\left(\frac{\|Z\|_{\Phi,\xi}}{\|\|Z\|_{\Phi,\xi}\|_\Phi}\right),
\end{eqnarray*}
We conclude by using the definition of the norm $\|\cdot\|_\Phi$.
\end{proof}
\begin{Lem}\label{phix}
If the Orlicz function $\Phi$ satisfies \eqref{submult} then  for any $E$-valued random variable $X$ we have $\|X\|_1\le \|X\|_\Phi$.
\end{Lem}
\begin{proof}
Using Jensen's inequality, we obtain
$$
\E \Phi\left(\frac{\|X\|}{\|X\|_1}\right)\ge \Phi(1).
$$
Note that $\Phi(1) \le\Phi(1)^2$ by \eqref{submult} and then that $\Phi(1)\ge 1$. We conclude that $\|X\|_1\le \|X\|_\Phi$ by using the definition of the norm $\|\cdot\|_\Phi$.
\end{proof}
\begin{Lem}\label{preleminary}
Let $u_0\ge 0$ and $(u_n)_{n\in \Z}$ be a real sequence such that
$|u_n|\le u_0$ if $n<0$. Assume that
\begin{equation}\label{recseq}
u_n= \sum_{i=1}^pa_iu_{n-i}, \qquad \forall n\ge0,
\end{equation}
where $a_1,\dots,a_p$ are fixed nonnegative numbers with
$a=\sum_{i=1}^p a_i<1$. Then, $$\sup_{k\ge n}u_k\le a^{n/p}u_0,
\qquad \forall n\ge0.$$
\end{Lem}
\begin{proof}
By a recursion argument, one first shows that $\sup_{k\le n} u_k\le u_0$. Then $(u_n)_{n\in\N}$ is bounded by $u_0$. Let $v_n=\sup_{k\ge n} u_k$ for $n\in\Z$. Using the relation (\ref{recseq}), we get $v_n\le a\, v_{n-p}$ for all $n\ge 0$. Then recursively $v_n\le a^{-[-n/p]}v_{n+p[-n/p]}$. From $|u_n|\le u_0$ if $n<0$, $v_{n+p[-n/p]}=v_0=u_0$ because $n+p[-n/p]\le 0$. The result follows from $-[-n/p]\ge n/p$.
\end{proof}
\subsection{$p$-Markov stationary approximations}
\label{approx}
In order to construct a solution to \eqref{chaineinf}
we consider, for each fixed $p\ge0$ and $q>0$, the $p$-Markov
process $(X_{p, q,t})_{t\ge0}$ defined by $X_{p, q,t}=0$ for
$t\le-q$ and  the recurrence equation
\begin{equation}
\label{chaineinfp}
 X_{p,q,t}=F(X_{p,q,t-1},\ldots,X_{p,q,t-p},0,0,\ldots;\xi_t) \qquad\mbox{if}\qquad t>q.
\end{equation}
Using the notation of lemma \ref{normxi} with $\xi=\xi_0$ and $\zeta=(X_{p,q,-1},X_{p,q,-2},\ldots)$ and $z(\xi,\zeta)=F(\zeta,\xi)$, the Lipschitz condition \eqref{lipcondf} implies that
$$\|X_{p,q+1,0}-X_{p,q,0}\|_{\Phi,\xi}\le \sum_{i=1}^pa_i\|X_{p,q+1,-i}-X_{p,q,-i}\|.$$
Applying lemma \ref{preliminary2},
\begin{eqnarray*}
\left\|X_{p,q+1,0}-X_{p,q,0}\right\|_\Phi&\le&\left\|\|X_{p,q+1,0}-X_{p,q,0}\right\|_{\Phi,\xi}\|_\Phi\\
&\le&\left\|\sum_{i=1}^pa_i\|X_{p,q+1,-i}-X_{p,q,-i}\|\right\|_\Phi\\
&\le&\sum_{i=1}^pa_i\left\|X_{p,q+1,-i}-X_{p,q,-i}\right\|_\Phi\\
&\le&\sum_{i=1}^pa_i\left\|X_{p,q+1-i,0}-X_{p,q-i,0}\right\|_\Phi.
\end{eqnarray*}
The last inequality follows from the fact that by the definition
of $X_{p,q,-i}$ and $X_{p,q-i,0}$, these quantities have the same law for each triplet of
positive integers $(p,q,i)$.
We now consider $v_n=\left\|X_{p,n+1,0}-X_{p,n,0}\right\|_\Phi$ for $n\in\Z$, with
$v_n=0$ if $n<0$. For $n>0$
$$
v_n\ \le \ \sum_{i=1}^pa_iv_{n-i}.
$$
From lemma \ref{preleminary} we obtain
$$
v_n\le a^{n/p}v_0\le a^{n/p}\|X_{p,1,0}\|_\Phi\le
a^{n/p}\|F(0,0,\ldots;\xi_t)\|_\Phi\le a^{n/p}\mu_\Phi.
$$
Hence, for each $p$, $(X_{p,n,0})_{n\in\N}$ is a Cauchy sequence in $\L^\Phi$; it converges to some $X_{p,0}\in \L^\Phi$. From its construction, it is clear that $X_{p,n,0}$ is measurable with respect to  the $\sigma-$algebra generated by $\{\xi_t, t\le0\}$.
The $\L^\Phi$-convergence ensures that this is also the case for $X_{p,0}$.
Hence there exists some
measurable function $H_p$ such that $X_{p,0} = H_p(\xi_0,\xi_{-1},\ldots)$. As $n\uparrow\infty$, a continuity argument on $F$ implies that
$X_{p,0}=F(X_{p,-1},\ldots,X_{p,-p},0,0,\ldots;\xi_0)$ and shifting the lag $t\in\Z$
leads to the equalities,
 $$X_{p,t} = H_p(\xi_t, \xi_{t-1},\xi_{t-2},\ldots)=F(X_{p,t-1},\ldots,X_{p,t-p},0,0,\ldots;\xi_t).$$
Then
the sequence $(X_{p,t})_{t\in\Z}$ is  a stationary solution of the
recurrence equation (\ref{chaineinfp}) for each $p\ge 0$.
\\
Consider
$$\mu_{\Phi,p}= \|X_{p,t}\|_\Phi ,
\qquad
\Delta_{\Phi,p,t}=\|X_{p+1,t}-X_{p,t}\|_\Phi ,
$$
The definition of $\mu_{\Phi,p}$ given here for $p>0$  extends to $p=0$ since $X_{0,t}=F(0,0,\dots;\xi_t)$ satisfies $\|X_{0,t}\|_\Phi=\mu_\Phi$ by \refeqn{propmom}.
\begin{Lem}
\label{bornemuq}\label{borneDeltaq}
Assume conditions \eqref{propinf} and \eqref{propmom} hold for some
Orlicz function $\Phi$ satisfying \eqref{submult}. Then
$$  \mu_{\Phi,\infty}=\sup_{p\ge0}\mu_{\Phi,p}\le \frac{ \mu_\Phi}{1-a}\qquad\mbox{and}\qquad  \Delta_{\Phi,p}=\sup_{t\in\Z}\Delta_{\Phi,p,t}\le a_{p+1} \frac{\mu_\Phi}{(1-a)^2}.$$
\end{Lem}
\begin{proof}
From \refeqn{propinf}, we have that
\begin{eqnarray*}
\mu_{\Phi,p}\le \|X_{p,t}-X_{0, t}\|_\Phi+\mu_\Phi\le\sum_{j=1}^{p}a_j\|X_{p,t-j}\|_\Phi+\mu_\Phi
\le\mu_{\Phi,p}\sum_{j=1}^{p}a_j+\mu_\Phi, \end{eqnarray*}
hence $\mu_{\Phi,p}\le(1-a)^{-1}\mu_\Phi$ and $\mu_{\Phi,\infty}\le(1-a)^{-1}\mu_\Phi$ follow.
In a similar way, we obtain the inequalities
\begin{eqnarray*}
   \Delta_{\Phi,p,t}&=& \Big\|F(X_{p+1,t-1},\ldots,X_{p+1,t-p-1},0,0,\dots;\xi_t)-\ F(X_{p,t-1},\ldots,X_{p,t-p},0,0,\dots;\xi_t)\Big\|_\Phi\\
   &\le&
   \sum_{j=1}^{p}a_j\|X_{p+1,t-j}-X_{p,t-j}\|_\Phi+a_{p+1}\|X_{p+1,t-p-1}\|_\Phi
   \\
&\le&
   \sum_{j=1}^{p}a_j\Delta_{\Phi,p, t-j}+a_{p+1}\|X_{p+1,0}\|_\Phi.
\end{eqnarray*}
This implies that $\Delta_{\Phi,p}\le a_{p+1}(1-a)^{-1} \mu_{\Phi,p+1}$ and the result of lemma \ref{bornemuq} is shown.
\end{proof}

\subsection{Proof of the existence of a solution to \eqref{chaineinf}}\label{sec::chaineinf}

Note first that lemma \ref{borneDeltaq} implies that $X_{p,t}\to_{p\to\infty}
X_{t}$ in $\L^\Phi$ since this space is complete. The continuity of $F$
ensures that $X_{t}$ is a solution of \refeqn{chaineinf}.
Furthermore, as a limit in $\L^\Phi$ of strictly stationary processes, $X_t$ is also stationary (in law) and
$\|X_t\|_\Phi<\infty$. Finally, $X_t = H(\xi_t,
\xi_{t-1},\ldots)$ is the limit in $\L^\Phi$ of $X_{p,t} =
H_p(\xi_t, \xi_{t-1},\ldots)$.

\subsection{Proof of the weak dependence properties}\label{wd::chaineinf}

The weak dependence property of a solution to \eqref{chaineinf} is formulated in terms of the $\L^1$-norm in the definition of the coefficients $\tau$. As shown in lemma \ref{phix},  $\|X\|_1\le \|X\|_\Phi$ for any $E$-valued random variable $X$. Then assumptions \eqref{propinf} and \eqref{propmom} are always satisfied replacing $\|\cdot\|_\Phi$ with $\|\cdot\|_1$. We first prove corollary \ref{The finite memory case}:
\begin{proof} We use coupling techniques
to evaluate the coefficients $\tau$, see p.\pageref{coupling}. Let $(\xi'_t)_{t\in\Z}$ be an independent copy of $(\xi_t)_{t\in\Z}$. We define the process $(X^\ast_{p,t})_{t\in\Z}$ as
$$X^\ast_{p,t}=\left\{
\begin{array}{cc}
F(X^\ast_{p,t-1},\ldots,X^\ast_{p,t-p},0,0,\ldots;\xi_t'),&\text{ for
}t\le0;\\
F(X^\ast_{p,t-1},\ldots,X^\ast_{p,t-p},0,0,\ldots;\xi_t), &\text{ for
}t>0;
\end{array}\right.
$$
Using similar arguments as section \ref{approx}, there exists a sequence of measurable variables with respect to  the $\sigma-$algebra generated by $\xi'_t, t\le0$ denoted by $(X^\ast_{p,n,0})_{n\in\N}$ such that it converges in $\L^\Phi$ to $X^\ast_{p,0}\in \L^\Phi$.
The $\L^\Phi$-convergence ensures that $X^\ast_{p,0}$ are also measurable variables with respect to  the $\sigma-$algebra generated by $\xi'_t, t\le0$.
Then, by definition of $\xi'_t, t\le0$, $X^\ast_{p,0}$ is independent of $X_{p,0}$. If there exists a non-increasing function $\delta_p(r)$ of $r$ such that $\|X_{p,r}-X_{p,r}^\ast\|_1\le \delta_p(r)$, we have $\tau_{p,r}\le\delta_p(r).$ This follows from the coupling property of weak dependence coefficients $\tau$ explained in \cite{DEPR04}.

\noindent
Assumption (\ref{propinf}) and lemma \ref{phix} yield
$$
\|X_{p,r}-X_{p,r}^\ast\|_1\le
\sum_{i=1}^pa_i\|X_{p,r-i}-X_{p,r-i}^\ast\|_1.
$$
Denoting $w_r=\|X_{p,r}-X_{p,r}^\ast\|_1$ for $r\in\Z$, we again use
lemma \ref{preleminary} and the relation  $\|F(0,0,\ldots;\xi_0)\|_1=\mu_1$ to obtain
$$w_r
\le a^{r/p}w_0
\le 2\mu_1 a^{r/p}
\le 2\frac{\mu_1}{1-a} a^{r/p}.$$
Now choose $\delta_p(r):=2\mu_1(1-a)^{-1} a^{r/p}$ leads to the result of corollary \ref{The finite memory case}.
\end{proof}
\noindent Now we finish the proof of theorem \ref{thchaineinf} defining the process $(X_t^\ast)_{t\in\Z}$ as the solution of the equations
$$X^\ast_{t\in\Z}=\left\{
\begin{array}{cc}
F(X^\ast_{t-1},X^\ast_{,t-2},\ldots;\xi_t'),&\text{ for
}t\le0;\\
F(X^\ast_{t-1},X^\ast_{t-2},\ldots;\xi_t), &\text{ for }t>0;
\end{array}\right.$$
We remark that $(X^\ast_t)_t$ is also a stationary chain with
infinite memory. Lem\-ma \ref{borneDeltaq} gives
$$\|X_r-X_{p,r}\|_1\le \sum_{k=p}^\infty  \Delta_{1,k}\le \frac{\mu_1}{(1-a)^2}\sum_{k=p}^\infty
a_{k+1}.
$$
The same bound holds for the quantity $\|X_r^\ast-X_{p,r}^\ast\|_1$.
For each integer $p$,
\begin{eqnarray*}\|X_r-X_r^\ast\|_1
\le
\|X_r-X_{p,r}\|_1+\|X_{p,r}-X^\ast_{p,r}\|_1+\|X_{r}^\ast-X^\ast_{p,r}\|_1
\le 2\frac{\mu_1}{1-a}\left(a^{r/p}+\sum_{k=p+1}^\infty  \frac{a_{k}}{1-a}\right).
\end{eqnarray*}
Because this bound is non-increasing with $r$, we conclude the weak dependence properties in theorem \ref{thchaineinf} by using the coupling technique.

\subsection{Proof of theorem \ref{clt}}\label{pfclt}
First we recall the assumption \eqref{Dq} of \cite{DEDO03} for $q>1$,
\begin{equation}\label{Dq}
\tag{D(q)}\int_0^{\|X_0\|_1}((\tau/2)^{-1}(u))^{q-1}Q^{q-1}\circ G(u)du<\infty,
\end{equation}
where $(\tau/2)^{-1}(u)=\inf\{k\in \N/\ \tau(k)\le 2u\}$. Here $Q$ denotes the generalized inverse of the tail function $x\mapsto\P(|X_0|>x)$ and $G$ the inverse of $x\mapsto \int_0^xQ(u)du$. Dedecker and Doukhan proved in \cite{DEDO03} the SLLN and the CLT under \eqref{Dq} for respectively $1<q<2$ and $q=2$. The SIP is proved in \cite{DEPR04} under \eqref{Dq} for $q=2$. Write $A(p)=\sum_{j\ge p} a_j$ and $A^{-1}$ its generalized inverse $A^{-1}(u)=\inf\{k\in\N/\ A(u)\le u\}$,
$$\Psi_q(x)=\Phi\big(x^{1/(q-1)}\big)/x^{1/(q-1)}\mbox{ and }
\displaystyle\Psi_q^\ast(x)=\sup_{y\ge 0}\{xy-\Psi_q(y)\}.
$$
Noticing that $A^{-1}(u)=k$ on $]A(k-1);A(k)]$ and that $  \widetilde\Phi_q(x)=\Psi_q^\ast(x^{q-1})$, there exists $C>0$ such that
$$
\int_0^a  \widetilde\Phi_q\left(c_0(A^{-1}(u)-1)\ln(u)\right)du\le C \sum_{k\ge1}a_k   \widetilde\Phi_q
\left(c_0k
\Big(1-\1_{ \{\sum_{j\ge k}a_j>0 \}}\ln\Big(\sum_{j\ge k}a_j\Big)\Big)\right).
$$
Then assumption \eqref{Dp} implies that we work under the condition
\begin{equation}\label{Dq'}
\int_0^a\Psi_q^\ast\left(\left(c_0(A^{-1}(u)-1)\ln(u)\right)^{q-1}\right)du<\infty.
\end{equation}
We want to prove that condition \eqref{Dq'} implies \eqref{Dq} for all $q> 1$. The first step is to prove the bound
\begin{equation}\label{boundtau-1}
(\tau/2)^{-1}(u)\le \left[\left(A^{-1}\left(\frac{1-a}{2\mu_1}u\right)-1\right)\frac{\ln \left(\frac{1-a}{2\mu_1}u\right)}{\ln a}\right],
\end{equation}
Theorem \ref{thchaineinf} gives
$
(\tau/2)^{-1}(u)\le \inf B$ with $$B=\left\{k\in\N\mbox{ such that }\exists p\ge1\mbox{ with } \frac{\mu_1}{1-a}\left(a^{k/p}+A(p+1)\right)\le u\right\}.$$
Set $v=(1-a)(2\mu_1)^{-1}u$, the integer $p^\ast=A^{-1}(v)-1$ is close to the infimum of $B$. Then all integers $k$ with $a^{k/p^\ast}\le v$ belong to $B$, as for instance $k^\ast=[(A^{-1}(v)-1)\ln v/\ln a]$ which is then larger than $
(\tau/2)^{-1}(u)$ by definition. Observe that $A^{-1}(v)=1$ as soon as $v\ge a$, thus $[(A^{-1}(v)-1)\ln v/\ln a]=0$ for $v\ge a$.\\

\noindent Using this estimate of $(\tau/2)^{-1}$ in \eqref{boundtau-1}, condition \eqref{Dq} holds if
\begin{equation}\label{integral}
\int_0^{a}\left[\frac{(A^{-1}(v)-1)\ln v}{\ln a}\right]^{p-1}Q^{p-1}\circ  G\left(\frac{2\mu_1}{1-a}v\right)dv<\infty.
\end{equation}
Let $\widetilde \Psi$ be an Orlicz function and $\widetilde \Psi^\ast(x)=\sup_{y>0}\{xy-\widetilde\Psi(y)\}$ be its Young dual function.
For any functions $f$ and $g$, Young's inequality gives:
\begin{multline*}
\int_0^af(x)g(x)dx\le 2\inf\left\{c>0\mbox{ with }\int_0^a\widetilde\Psi\left(\frac{f(x)}{c}\right)dx\le 1\right\}\\
\times\inf\left\{c>0\mbox{ with }~\int_0^a\widetilde\Psi^\ast \left(\frac{g(x)}{c}\right)dx\le 1\right\}.
\end{multline*}
In the following we apply this inequality with and $\widetilde \Psi=K\Psi_p$ for some $K>0$, $f(x)=Q^{p-1}\circ G(2\mu_1 (1-a)^{-1}x)$ and $g(x)\le ((A^{-1}(x)-1)\ln(1/x)(-\ln a)^{-1})^{q-1}$. Note that the Young dual function is here $\widetilde \Psi^\ast(x)=K\Psi_q^\ast(x/K)$ and then $\int_0^af(x)g(x)dx$ is equal to the left hand side term \eqref{integral} up to the choice of the constant $K>0$, see below.
In view of Young's inequality, the first term in the bound of \eqref{integral} thus expresses as the infimum over $c>0$ such that
$$
K\frac{1-a}{2\mu_1}\int_0^{\|X_0\|_1}\frac{\Phi\left(Q\circ G\left(u\right)/c\right)}{Q\circ G\left(u\right)/c}du\le1.
$$
Replacing $G(u)$ with $x$, one obtains the simpler inequality:
$$
K\frac{1-a}{2\mu_1}\int_0^1\Phi\left(\frac{Q(x)}{c}\right)cdx=K\frac{1-a}{2\mu_1}c\E\Phi\left(\frac{|X_0|}{c}\right)\le 1.
$$
The last equality is set using the definition of $Q(x)$. If assumption \eqref{propmom} holds, the last inequality is satisfied for $K= 2\mu_1\mu_\Phi^{-1}$ and $c=\mu_\Phi(1-a)^{-1}$.\\

\noindent The second term of the Young inequality expresses as the infimum over $c>0$ such that
\begin{equation}\label{cond}
K\int_0^a\Psi^\ast_q\left(\frac{((A^{-1}(x)-1)\ln(1/x))^{q-1}}{K(-\ln a)^{q-1}c}\right)dx\le 1.
\end{equation}
Because $  \widetilde\Phi_q(x)=\Psi_q^\ast(x^{q-1})$ we check that
$$\displaystyle
0<\frac{\int_0^a  \widetilde\Phi_q\left(c_0(A^{-1}(u)-1)\ln(u)\right)du\vee1}{(K\wedge 1)(-\ln a)^{p-1}}=:c_1$$ satisfy the relation \eqref{cond}. It is obvious by \eqref{Dq'} that $c_1<\infty$ and then we proved the implications
$$\eqref{Dp}\mbox{ with }q>1\Rightarrow \eqref{Dq'}\mbox{ with }q>1\Rightarrow \eqref{Dq}.
$$
This ends the proof as the results of theorem \ref{clt} are versions of the results in \cite{DEDO03,DEPR04} that hold under assumption \eqref{Dq}.

\subsection{Proof of proposition \ref{simulation}}\label{pfsimu}
Let $n$ be a fixed integer and $s_n\le n-1$. Let $(X_t)_{t\in\Z}$ be the stationary solution of $X_{t}=F(X_{t-1},X_{t-2},0,0,\ldots;\xi_t)$. The Lipschitz assumption \eqref{lipcondf} implies for $1\le k\le n$
$$
\left\|\widetilde
X_k-X_{k}\right\|_\Phi\le\sum_{i=1}^{k-1}a_i\left\|\widetilde
X_{k-i}-X_{k-i}\right\|_\Phi+\sum_{i\ge k}a_i\|X_0-c_i\|_\Phi.
$$
The sequence $v_{k}=\left\|\widetilde
X_{k+1}-X_{k+1}\right\|_\Phi$, $k=1,2,\ldots$
satisfies the recursion
$$v_k\le \sum_{j=1}^ka_j v_{k-j}+u_k \qquad \mbox{ for all }\qquad k\ge1
$$
where $u_k=(\|X_0\|_\Phi+\overline c)\sum_{j>k}a_j$ for $k\ge1$. Notice that $u_k\downarrow_{k\to\infty}0$. We first prove the boundedness of $(v_k)_{k\in\N}$. Let $\ell$ be a fixed integer. For all $k$ such that $\ell \ge k$, $v_k\le a \sup_{i\le \ell}v_i+u_1$. We deduce that $\sup_{i\le \ell}v_i\le u_1$. Finally $\|v\|_\infty\le a\|X_0\|_\Phi/(1-a)$.
\\
Now for all integers $k,s \ge 1$ such that $\ell \ge k+s$,
\begin{eqnarray*}
v_\ell\le \sum_{j=1}^k a_jv_{\ell-j}+\sum_{j=k+1}^\ell a_jv_{\ell-j}+u_\ell
\le a\sup_{j\ge s}v_j+\|v\|_\infty\sum_{j=k+1}^\infty a_j+u_{k+s}.
\end{eqnarray*}
This inequality holds for all $\ell \ge k+s$. Then
$$
\sup_{j\ge k+s}v_j\le a\sup_{j\ge s}v_j+\|v\|_\infty\sum_{j=k+1}^\infty a_j+u_{k}.
$$
We deduce that
$$
\sup_{j\ge nk}v_j\le a^n\|v\|_\infty+\frac{1}{1-a}\left(\|v\|_\infty\sum_{j=k+1}^\infty a_j+u_{k}\right).
$$
Using the inequality $\|v\|_\infty\le a\|X_0\|_\Phi/(1-a)$, one gets the result.

  \begin{center}
    {\bf APPENDIX}
  \end{center}
We give below general conditions for the existence and the boundedness of joint densities of Affine Models defined in section \ref{secexample}. Thus we extend the results for Bilinear Models given in \cite{DOMA05}.
\begin{Prop}[Regularity of affine models]
\label{propregaf}
\normalfont

Here $E=E'=\R^d$ for some $d\ge 1$. Suppose that the innovations $(\xi_t)_{t\in\Z}$  in the model
(\ref{modaf}) have a common bounded marginal density $f_\xi$. Moreover, if $\inf_{(x_j)_{j>0}}\det M((x_j)_{j>0})=\underline M>0$, the marginal densities
$f_{X_1,\dots,X_n}$ of $(X_1,\dots,X_n)$ exist for all $n>0$
and satisfy
$$\|f_{X_1,\dots,X_n}\|_\infty\le \underline M ^{-n}\|f_\xi\|_\infty^n.$$
\end{Prop}
\begin{proof}
The solution $X_t = H(\xi_t,
\xi_{t-1},\ldots)$ obtained in section \ref{sec::chaineinf} is independent of $(\xi_j)_{j>t}$. If $G_1$ is a bounded continuous function on $E$ with value in $\R$, it holds that
\begin{eqnarray*}
\E \,G(X_1)&=&\E\, G_1(M(X_{0},\ldots)\xi_1+f(X_{0},X_{-1},\ldots))\\
&=&\int\int \,G_1(M(u)s_1+f(u))\,f_\xi (s_1)ds_1\P_{(X_{0},X_{-1},\ldots)}(du)\\
&\le&\underline M\int\int \,G(x_1)\,f_\xi (M^{-1}(u)(x_1-f(u)))\P_{(X_{0},X_{-1},\ldots)}(du)ds_1.
\end{eqnarray*}
The last inequality follows by the substitution $M(u)s_1+f(u)=x_1$ valid under the assumption $\inf_{(x_j)_{j>0}}\det M((x_j)_{j>0})=\underline M>0$ ensuring that $M(u)$ is invertible for all $u$. We obtain
$$
f_{X_1}(x_1)\le\underline M^{-1} \int f_\xi (M^{-1}(u)(x_1-f(u)))\P_{(X_{0},X_{-1},\ldots)}(du)\le \underline M^{-1} \|f_\xi\|_\infty.
$$
We proceed by induction for the cases $n\ge 2$. Assume that $\|f_{X_1,\dots,X_{n-1}}\|_\infty\le \underline M
^{-(n-1)}\|f_\xi\|_\infty^{n-1}$ is satisfied. Let $G_n$ be a bounded continuous function on $E^n$ with value in $\R$, one has
\begin{eqnarray*}
\E \,G_n(X_1,\dots,X_n)&=&\E\, G_n(X_1,\dots,X_{n-1},M(X_{n-1},X_{n-2},\ldots)\xi_n+f(X_{n-1},X_{n-2},\ldots))\\
&&\hspace{-4cm}=\int\int\int
G_n(x_1,\dots,x_{n-1},M(x_{n-1},\dots,x_1,u)s_n
+f(x_{n-1},\dots,x_1,u))\\&&\hspace{-4cm}f_\xi (s_n)ds_nf_{(X_1,\dots,X_{n-1})}(x_1,\dots,x_{n-1})dx_1\cdots dx_{n-1}dP_{(X_{0},X_{-1},\ldots|X_1,\dots,X_{n-1})}(u).
\end{eqnarray*}
The substitution $M(x_{n-1},\dots,x_1,u)s_n
+f(x_{n-1},\dots,x_1,u)=x_n$ yields
\begin{multline*}
f_{X_1,\dots,X_n}(x_1,\dots,x_n)\le\hspace{-1mm}\underline M^{-1}\hspace{-2mm} \int\hspace{-3mm}\int f_\xi (M^{-1}(x_{n-1},\dots,x_1,u)(x_n-f(x_{n-1},\dots,x_1,u)))\\
f_{(X_1,\dots,X_{n-1})}(x_1,\dots,x_{n-1})dx_1\cdots dx_{n-1}dP_{(X_{0},X_{-1},\ldots|X_1,\dots,X_{n-1})}(u).
\end{multline*}
Together with the induction assumption $\|f_{X_1,\dots,X_{n-1}}\|_\infty\le \underline M
^{-(n-1)}\|f_\xi\|_\infty^{n-1}$, this last inequality yields $\|f_{X_1,\dots,X_{n}}\|_\infty\le \underline M
^{-n}\|f_\xi\|_\infty^{n}$.
\end{proof}

\noindent {\bf Acknowledgements.} We are deeply grateful to Alain
Latour, and Thomas Mikosch who made a critical review of the drafts
and with whom we have worked on the final version of this paper. We
also wish to thank J\'{e}r\^{o}me Dedecker for important comments relating
$\tau$ dependence conditions as well as Lionel Truquet who pointed
out the example of Galton-Watson processes with immigration. We
finally wish to thank an anonymous referee whose remarks considerably helped us
to improve our results.
\bibliographystyle{plain}
\bibliography{ChaineDW21}

\begin{thebibliography}{10}

\bibitem{ANDR84}
D.~W.~K. Andrews.
\newblock Nonstrong mixing autoregressive processes.
\newblock {\em J. Appl. Probab.}, 21(4):930--934, 1984.

\bibitem{BERB87}
H.~Berbee.
\newblock Chains with infinite connections: uniqueness and {M}arkov
  representation.
\newblock {\em Probab. Theory Related Fields}, 76(2):243--253, 1987.

\bibitem{BOUG93}
P.~Bougerol.
\newblock Kalman filtering with random coefficients and contractions.
\newblock {\em Probab. Theory Related Fields}, 31:942--959, 1993.

\bibitem{BUWY99}
P.~B{\"u}hlmann and A.~J. Wyner.
\newblock Variable length {M}arkov chains.
\newblock {\em Ann. Statist.}, 27(2):480--513, 1999.

\bibitem{COFE02}
F.~Comets, R.~Fern{\'a}ndez, and P.~A. Ferrari.
\newblock Processes with long memory: regenerative construction and perfect
  simulation.
\newblock {\em Ann. Appl. Probab.}, 12(3):921--943, 2002.

\bibitem{DEDO03}
J.~Dedecker and P.~Doukhan.
\newblock A new covariance inequality and applications.
\newblock {\em Stochastic Process. Appl.}, 106(1):63--80, 2003.

\bibitem{DEDO05}
J.~Dedecker, P.~Doukhan, G.~Lang, J.~R. Le\'on, S.~Louhichi, and C.~Prieur.
\newblock {\em Weak Dependence, Examples and Applications}, volume 190 of {\em
  Lecture Notes in Statistics}.
\newblock Springer-Verlag, Berlin, 2007.

\bibitem{DEPR04}
J.~Dedecker and C.~Prieur.
\newblock Coupling for $\tau$-dependent sequences and applications.
\newblock {\em J. Theor. Probab.}, 17(4):861--855, 2004.

\bibitem{DIFR99}
P.~Diaconis and D.~Freedman.
\newblock Iterated random functions.
\newblock {\em SIAM Rev.}, 41(1):45--76, 1999.

\bibitem{DOBR70}
R.~L. Dobrushin.
\newblock Prescribing a system of random variables by conditional
  distributions.
\newblock {\em Theory Probab. Appl.}, 15:458--486, 1970.

\bibitem{DOKU93}
R.~L. Dobrushin and S.~Kusuoka.
\newblock {\em Statistical Mechanics and Fractals}, volume 1567 of {\em Lecture
  Notes in Mathematics}.
\newblock Springer-Verlag, New York, 1993.

\bibitem{DOUK94}
P.~Doukhan.
\newblock {\em Mixing}, volume~85 of {\em Lecture Notes in Statistics}.
\newblock Springer-Verlag, New York, 1994.

\bibitem{DOLO99}
P.~Doukhan and S.~Louhichi.
\newblock A new weak dependence condition and applications to moment
  inequalities.
\newblock {\em Stochastic Process. Appl.}, 84(2):313--342, 1999.

\bibitem{DOMA05}
P.~Doukhan, H.~Madre, and M.~Rosenbaum.
\newblock {ARCH} type bilinear weakly dependent models.
\newblock {\em Statistics}, 41(1):31--45, 2007.

\bibitem{dmr}
P.~Doukhan, P.~Massart, and E.~Rio.
\newblock The functional central limit theorem for strongly mixing processes.
\newblock {\em Ann. Inst. H. Poincar{\'e} Probab. Statist.}, 30(2):63--82,
  1994.

\bibitem{DOTE06}
P.~Doukhan, G.~Teyssi{\`e}re, and P.~Winant.
\newblock A {LARCH}($\infty$) vector valued process.
\newblock In Patrice Bertail, Paul Doukhan, and Philippe Soulier, editors, {\em
  Dependence in Probability and Statistics}, volume 187 of {\em Lectures Notes
  in Statistics}. Springer, New York, 2006.

\bibitem{DUFL97}
M.~Duflo.
\newblock {\em Random Iterative Models}, volume~34 of {\em Applications of
  Mathematics}.
\newblock Springer-Verlag, Berlin, 1997.

\bibitem{FEMA04}
R.~Fern{\'a}ndez and G.~Maillard.
\newblock Chains with complete connections and one-dimensional {G}ibbs
  measures.
\newblock {\em Electron. J. Probab.}, 9(6):145--176, 2004.

\bibitem{GILE04}
L.~Giraitis, R.~Leipus, P.~M. Robinson, and D.~Surgailis.
\newblock \mbox{LARCH}, leverage and long memory.
\newblock {\em J. Financial Econometrics}, 2(2):177--210, 2004.

\bibitem{IOGR90}
M.~Iosifescu and S.~Grigorescu.
\newblock {\em Dependence with Complete Connections and its Applications},
  volume~96 of {\em Cambridge Tracts in Mathematics}.
\newblock Cambridge University Press, Cambridge, 1990.

\bibitem{IOTH69}
M.~Iosifescu and R.~Theodorescu.
\newblock {\em Random Processes and Learning}.
\newblock Springer-Verlag, New York, 1969.

\bibitem{KACM59}
M.~Kac.
\newblock {\em Probability and Related Topics in Physical Sciences}, volume~1a
  of {\em Lectures in Applied Mathematics Series}.
\newblock American Mathematical Society, London, 1959.

\bibitem{KALL97}
O.~Kallenberg.
\newblock {\em Foundations of Modern Probability}.
\newblock Probability and its Applications. Springer-Verlag, New York, 1997.

\bibitem{KR}
M.~A. Krasnoselskii and Y.~B. Rutickii.
\newblock {\em Convex Functions and Orlicz Spaces}.
\newblock Noordhoff Ltd., Groningen, 1961.

\bibitem{LALI80}
L.~D. Landau and E.~M. Lifshitz.
\newblock {\em Statistical Physics}, volume~5 of {\em Course of Theoretical
  Physics}.
\newblock Butterworth-Heinemann, Oxford, 3 edition, 1980.

\bibitem{Latour1997}
A.~Latour.
\newblock The multivariate \mbox{GINAR(p)} process.
\newblock {\em Adv. Appl. Prob.}, 29:228 -- 247, 1997.

\bibitem{Merlevede2006}
F.~Merlevède and M.~Peligrad.
\newblock On the weak invariance principle for stationary sequences under
  projective criteria.
\newblock {\em J. Theoret. Probab.}, 19:647--689, 2006.

\bibitem{Peligrad2005}
M.~Peligrad and S.~Utev.
\newblock A new maximal inequality and invariance principle for stationary
  sequences.
\newblock {\em Ann. Probab.}, 33:798--815, 2005.

\bibitem{RAWI06}
N.~Ragache and O.~Wintenberger.
\newblock Convergence rates for density estimators of weakly dependent time
  series.
\newblock In P.~Bertail, P.~Doukhan, and P.~Soulier, editors, {\em Dependence
  in Probability and Statistics}, volume 187 of {\em Lectures Notes in
  Statistics}. Springer, New York, 2006.

\bibitem{Ri}
E.~Rio.
\newblock {\em Th{\'e}orie asymptotique des processus al{\'e}atoires faiblement
  d{\'e}pendants}, volume~31 of {\em Math{\'e}matiques \& Applications}.
\newblock Springer-Verlag, Berlin, 2000.

\bibitem{Wu2005}
W.~B. Wu.
\newblock Nonlinear system theory: Another look at dependence.
\newblock {\em Proc. Natl. Acad. Sci. USA}, 102:14150--14154, 2005.

\end{thebibliography}
\end{document}